\documentclass[12pt]{article}
\usepackage{amsmath}
\usepackage{amscd}
\usepackage{graphicx}
\usepackage{amssymb,amsfonts,amsthm,amscd,mathrsfs}
\usepackage{verbatim}
\usepackage[T2A]{fontenc}
\usepackage[utf8]{inputenc}
\usepackage[russian]{babel}
\usepackage{bm}
\usepackage{authblk}

\renewcommand{\b}{\beta}

\newcommand{\e}{\varepsilon}

\newcommand{\E}{\mathsf{E}}

\begin{document}

\title{On the Power of Symmetrized Pearson's Type Test under Local Alternatives in Autoregression with Outliers}
\author{M.~V.~Boldin    
\footnote{Moscow State Lomonosov Univ., Dept. of Mech. and Math., Moscow, Russia\\
e-mail: boldin$_{-}$m@hotmail.com}}

\date{ }
\maketitle

\textbf{Abstract}

We consider a stationary linear AR($p$) model with observations subject to
gross errors  (outliers). The autoregression parameters are unknown as well as the distribution function $G$ of innovations. 
The distribution of outliers $\Pi$ is unknown and arbitrary, their intensity is
 $\gamma n^{-1/2}$ with an unknown $\gamma$, $n$ is the sample size.
 We test the hypothesis $H_0\colon G=G_0$ with simmetric $G_0$. We find the power of the test under local alternatives 
 $H_{1n}(\rho)\colon G=(1-\rho n^{-1/2})G_0+\rho n^{-1/2}H$. Our test is the special symmetrized Pearson's type test.  Namely, first of all we estimate the autoregression parameters  
  and then using the residuals from
the estimated autoregression we construct a kind of empirical distribution function
(e.d.f.), which is a counterpart of the (inaccessible)
e.d.f. of the autoregression innovations. We obtain a stochastic expansion
of this e.d.f. and its  symmetrized variant under $H_{1n}(\rho)$ , which enables us to construct and investigate our  symmetrized test of Pearson's type for $H_0$.
We establish qualitative robustness of this test in terms of uniform equicontinuity
 of the limiting power  (as functions of $\gamma,\rho$ and $\Pi$) with respect to $\gamma$ in a neighborhood of $\gamma=0$.



\section{
Введение и постановка задачи}

В этой работе мы рассматриваем стационарную AR($p$) модель с ненулевым средним
\begin{equation}
v_t = \b_1 v_{t-1} + \dots + \b_p v_{t-p} +\nu + {\e}_t, \quad  t \in \mathbb{Z} .
\end{equation}

В (1.1) $\{{\e}_t\}$ -- независимые одинаково распределенные случайные величины (н.о.р.сл.в.) с неизвестной функцией распределения (ф.р.) $G(x)$;
 $\E \e_1 = 0$,  $0<\E {\e}_1^2 < \infty$; $\bm{\b} = (\b_1, \dots, \b_p)^T \in  \mathbb{R}^p $-- вектор неизвестных параметров, таких что корни соответствующего (1.1) характеристического уравнения по модулю меньше единицы;  $\nu$ -- неизвестное среднее,  $\nu \in \mathbb R^1$.\\
 Эти условия дальше всегда предполагаютя выполненными и особо не оговариваются.\\
 Мы рассматриваем модель (1.1) с выбросами в наблюдениях. А именно, предполагается, что наблюдаются величины
$$
y_t=v_t+z^{\gamma_n}_t {\xi}_t,\quad t  = 1-p, \dots, n, \eqno(1.2)
$$
где $v_{1-p}, \dots, v_n$ -- выборка из стационарного решения $\{v_t\}$ уравнения  (1.1); $\{z^{\gamma_n}_t\}$ н.о.р.сл.в., принимающие значения 1 и 0, причем вероятность единицы
 $\gamma_n$,
 $$
 \gamma_n = \min(1, \frac{\gamma}{\sqrt{n}}) ,\quad  \gamma \ge 0 \text{\; неизвестно.}
 $$
Кроме того, $\{ \xi_t \}$ -- н.о.р.сл.в. с произвольным и неизвестным распределением $\Pi$.
Переменные $\{ \xi_t \}$ интерпретируются как выбросы (засорения), $\gamma_n$ уровень засорения. Для $\gamma = 0$ мы получаем модель (1.1) без засорений.

Модель (1.2) -- локальный вариант хорошо известной модели засорения данных во временных рядах, см. \cite{MartYoh86}.\\

 В настоящей работе мы рассматриваем гипотезу
 $$
 H_0\colon G(x)=G_0(x)\, \text{с известной функцией распределения }G_0(x).
 $$
 В качестве альтернативы к $H_0$ берется предположение о том, что $\{\e_t\}$  в (1.1) -- н.о.р.сл.в. с ф.р.
 $$
 A_n(x)=(1-\rho_n)G_0(x)+\rho_nH(x),\eqno(1.3)
 $$
 $$
 H(x)-\text{ф.р.}, \rho_n=\min(1,\frac{\rho}{\sqrt{n}}), \rho\geq 0.
 $$ 
 Предположение (1.3) будем понимать как локальную альтернативу к $H_0$ и обозначать $H_{1n}(\rho)$.\\
 Мы  хотим   построить тест типа хи-квадрат Пирсона для проверки $H_0$ иисследовать асимптотическую мощность этого теста при $H_{1n}(\rho)$. В частности, мы хотим установить качественную робастность асимпттической локальной мощности.\\
 Для решения этой задачи мы строим оценку ф.р. $G(x)$, она называется остаточной эмпирической функцией распределения (о.э.ф.р.). Мы получаем стохастическое разложение о.э.ф.р. при $H_{1n}(\rho)$. Этот результат имеет самостоятельное значение и вот почему.\\
В \cite{Bold.Arx.2} был построен тест для проверки нормальности инноваций в схеме (1.1)--(1.2), т.е. гипотезы
$
H_{\Phi}\colon G(x) \in \{\Phi(x/\theta),\; \theta > 0\},
$
где $\Phi(x)$ -- стандартная нормальная ф.р. Напомним, кстати, что гипотеза $H_{\Phi}$ эквивалентна нормальности самой стационарной последовательности 
 $\{v_t\}$. Нормальность инноваций обеспечивает оптимальность многих процедур оценивания и проверки гипотез в авторегрессии, см., например, \cite{And.},\cite{Broc.Dav.}. Поэтому проверка $H_{\Phi}$ -- содержательная задача. 
 В \cite{Bold.Arx.2} была построена статистика типа хи-квадрат Пирсона для $H_{\Phi}$ в схеме (1.1) --(1.2) и найдено ее распределение при гипотезе. Необходимый следующий шаг -- найти мощность построенного теста при локальных альтернативах и исследовать устойчивость мощности к выбросам. 
 Настоящая работа -- шаг в этом напрвлении, окончательный результат будет опубликован в отдельной заметке.
 
  Отметим, что недавно в \cite{Bold.Petr.},\cite{Bold.2020} были получены результаты о проверке $H_0$ и  $H_{\Phi}$ в схеме (1.1)--(1.2)  при $\nu=0$. Нулевое среднее -- существенное предположение.
 В нашей ситуации при неизвестном и, быть может, ненулевом $\nu$, использовать статистики типа хи квадрат из \cite{Bold.Petr.}, \cite{Bold.2020} не удается, т.к. их предельное распределение зависит от оценок $\bm{\b}$ и $\nu$ даже при $\gamma=0$ .\\
 Мы преодолеваем эту трудность в случае симметричной $G_0(x)$, строя специальную симметризованную статистику Пирсона, она является функционалом от симметризованной о.э.ф.р. Асимптотические распределения симметризованных статистик при $H_{0}$ и  $\gamma=0$ свободны.  
 Кроме того, симметризованный тест Пирсона оказывается качественно робастным.
 
 Все определения и результаты (основные, это Теоремы 2.1--2.2) представлены в Разделе 2.
 
 \section{Основные результаты}

\subsection{Стохастическое разложение о.э.ф. р. }

Относительно альтернативы $H_{1n}(\rho)$ из (1.3) нам потребуется сделать два предположения. Напомним, при $H_{1n}(\rho)$ ф.р. инноваций есть смесь вида
 $$
 A_n(x)=(1-\rho_n)G_0(x)+\rho_nH(x),
 $$
 $$
 H(x)-\text{ф.р.}, \rho_n=\min(1,\frac{\rho}{\sqrt{n}}), \rho\geq 0.
 $$

{\bf Условие (i).} Случайные величины с функциями распределения $G_0(x)$ и $H(x)$ имеют нулевые средние и конечные дисперсии.

{\bf Условие (ii).} Функции распределения $G_0(x)$ и $H(x)$ дважды дифференцируемы с ограниченными вторыми производными.

Перепишем уравнение (1.1) в удобном для дальнейшего рассмотрения виде. Для этого определим константу $\mu$ соотношеним
$$
\nu=(1-\b_1-\ldots-\b_p)\mu,
$$
тогда
$$
v_t -\mu= \b_1 (v_{t-1}-\mu) + \dots + \b_p( v_{t-p} -\mu)+ {\e}_t, \quad  t \in \mathbb{Z.}
$$

Если положить $u_t:=v_t-\mu,$ то
$$
v_t=\mu+u_t,\quad u_t=\b_1 u_{t-1} + \dots + \b_pu_{t-p}+ {\e}_t, \quad  t \in \mathbb{Z.}\eqno(2.1	)
$$
Последовательность $\{u_t\}$ в (2.1) -- авторегрессионная последовательность с нулевым средним и конечной дсперсией.\\
Построим по наблюдениям $\{y_t\}$ из (1.2) оценки ненаблюдаемых $\{\e_t\}$.\\ 
Далее $\Gamma,\, R$ --любые конечные неотрицательные числа. Пусть $\hat\mu_n$ будет любая последовательность, для которой при $H_{1n}(\rho)$ последовательнсть
$$
n^{1/2}(\hat\mu_n-\mu)=O_P(1),\quad n\to\infty, \text{ равномерно по}\,\,\gamma \leq\Gamma, \,\rho\leq R.\eqno(2.2)
$$
 Положим
$$
\hat u_t=y_t-\hat\mu_n,\quad t=1,\ldots,n.
$$
Пусть $\hat{\bm\b}_n= (\hat \b_{1n}, \dots, \hat \b_{pn})^T$ будет любая последовательность,
для которой при $H_{1n}(\rho)$ последовательность\\
 $$
 n^{1/2}(\hat{\bm\b}_n-\bm\b)=O_P(1),\quad n\to \infty, \text{ равномерно по}\,\,\gamma\leq \Gamma,\,\rho\leq R.\eqno(2.3)
 $$
Примеров подходящих оценок $\hat\mu_n,\,\hat{\bm\b}_n$ много, широкий класс составляют, например, М-оценки, построенные по $\{y_t\}$ и $\{\hat u_t\}$ аналогично тому, как они строятся по незасоренным данным, см., например,   \cite{Koul87}. Результаты для М-оценок при альтернативе $H_{1n}(\rho)$ в схеме (1.1)--(1.2) вполне аналогичны изложенным в Разделе 2.4 \cite{Bold.Petr.}. \\
Положим
$$
\hat \e_t =  \hat {u}_t -  \hat \b_{1n}\hat { u}_{t-1} - \dots - \hat \b_{pn}\hat{ u}_{t-p},\quad t = 1, \dots, n.\eqno(2.4)
$$
Величины $\{\hat \e_t\}$ называются остатками. 
Наша ближайшая цель -- исследовать асимптотические свойства о.э.ф.р.
$$
\hat G_n(x) = n^{-1} \sum_{t=1}^n I(\hat \e_t \le x),\quad  x\in \mathbb{R}^1.
$$
 Здесь и в дальнейшем $I(\cdot)$ обозначает индикатор события.\\
Функция $\hat G_n(x)$ -- аналог гипотетической э.ф.р.
$$
G_n(x) = n^{-1} \sum_{t=1}^n I(\e_t \le x)
$$
ненаблюдаемых величин $\e_{1}, \dots, \e_n$.\\ Легко проверить, что
$$
\hat \e_t=\e_t-(\hat{\bm\b}_n-\bm\b)^T
\tilde {\bm u}_{t-1}-
\hat\delta_n(\hat\mu_n-\mu)
+\alpha_{tn}-\sum_{j=1}^p(\hat \beta_{jn}-\beta_j)z_{t-j}^{\gamma_n} \xi_{t-j}.\eqno (2.5)
$$
В (2.5) 
$$
\tilde {\bm u}_{t-1}:=(u_{t-1},\ldots,u_{t-p})^T,\,\alpha_{tn}:=z_{t}^{\gamma_n} \xi_{t}-\sum_{j=1}^p \beta_{j}z_{t-j}^{\gamma_n} \xi_{t-j},\,\hat\delta_n:=1-\hat\b_1-\ldots-\hat\b_p.
$$
В силу 2.5)
$$
\hat G_n(x) = n^{-1} \sum_{t=1}^n I( \e_t \le x+(\hat{\bm\b}_n-\bm\b)^T\tilde {\bm u}_{t-1}+\hat\delta_n(\hat\mu_n-\mu)-\alpha_{tn}+\sum_{j=1}^p(\hat \beta_{jn}-\beta_j)z_{t-j}^{\gamma_n} \xi_{t-j}).\eqno (2.6)
$$
Пусть $
\bm{ \tau}_1=(\tau_1,\ldots,\tau_p)^T \in \mathbb{R}^p,\,\tau_0 \in \mathbb{R}^1
$. Свяжем с $\hat G_n(x)$ из (2.6) функцию
$$
 G_n(x,\tau_0,\bm{ \tau}_1) = n^{-1} \sum_{t=1}^n I( \e_t \le x+n^{-1/2}\bm{ \tau}_1^T\tilde {\bm u}_{t-1}+n^{-1/2} \tau_0-\alpha_{tn}+\sum_{j=1}^p n^{-1/2}\tau_j z_{t-j}^{\gamma_n} \xi_{t-j}).
$$
Очевидно, при $\gamma=0$ $G_n(x,0,\bm 0)=G_n(x)$.\\
С помощью тех же аргументов, что при доказательстве Теоремы 2.1 в \cite{Bold.Petr.} устанавливается следующая Теорема 2.1. В ней и далее $|\cdot|$ означает Евклидову норму вектора или матрицы.
\newtheorem{Th}{Теорема}[section]
\begin{Th}
Предположим, что верна альтернатива $H_{1n}(\rho)$. Пусть выполнены Условия (i) -- (ii), и $g_0(x)$ есть производная $G_0(x)$. Пусть $\Gamma \ge 0,\,R\ge 0$ and $\Theta \ge 0$ будут любые конечные числа. 
Тогда для любого $\delta > 0$
$$
\sup_{\gamma \le \Gamma,\rho\leq R} \Prob(\sup_{|\tau_0|\leq \Theta,|\bm {\tau}_1|\leq \Theta}|n^{1/2} [G_n(x,\tau_0,\bm{ \tau}_1)- G_n(x)] -g_0(x)\tau_0- \gamma\Delta(x, \Pi) | > \delta) \to 0,\quad n \to \infty.
$$
Здесь сдвиг
$$
\Delta(x, \Pi) = \sum_{j=0}^p [\E G_0(x + \b_j \xi_1) - G_0(x)],\, \b_0 =-1.
$$
\end{Th}
Теорема 2.1 и соотношения (2.2) --(2.3) стандартным образом (см. доказательство Следствия 2.1 в \cite{Bold.Petr.}) влекут
\newtheorem{Corollary}{Следствие
}[section]
\begin{Corollary}
При условиях Теоремы 2.1 
$$
\sup_{\gamma \le \Gamma,\rho\leq R} \Prob(|n^{1/2}[\hat G_n(x) -G_n(x)]- g_0(x)\hat\delta_n n^{1/2}(\hat\mu_n-\mu)-  \gamma\Delta_0(x, \Pi) | > \delta) \to 0,
\quad n \to \infty.
$$
\end{Corollary}

В связи с проверкой гипотезы $H_0$ нас особо интересует ситуация, когда ф.р. $G_0(x)$ симметрична относительно нуля. В этом случае будем брать оценкой $G(x)$ симметризованную оценку
$$
\hat S_n(x):=\frac{\hat G_n(x)+1-\hat G_n(-x)}{2}.\eqno(2.7)
$$
Положим
$$
\Delta_S(x, \Pi):=\frac{\Delta_0(x, \Pi)-\Delta_0(-x, \Pi)}{2}.\eqno(2.8)
$$
Пусть
$$
 S_n(x):=\frac{ G_n(x)+1- G_n(-x)}{2}.\eqno(2.9)
$$

Следствие 2.1 прямо влечет
\begin{Corollary}
При условиях Теоремы 2.1 в случае симметричной относительно нуля $G_0(x)$
$$
\sup_{\gamma \le \Gamma,\rho\leq R} \Prob(|n^{1/2}[\hat S_n(x) -S_n(x)]-  \gamma\Delta_S(x, \Pi) | > \delta) \to 0,
\quad n \to \infty.
$$
\end{Corollary}

\subsection{Тест типа хи-квадрат Пирсона для $H_{0}$}

В этом разделе мы построим тест типа хи-квадрат Пирсона для гипотезы 
$$
H_{0}\colon  G(x) =G_0(x), \,G_0(x)\,\text{известная симметричная относительно нуля ф.р.}
$$
Мы найдем асимптотическую мощность этого теста при локальных альтернативах $H_{1n}(\rho)$.

Будем  также предполагать $H(x)$ симметричной относительно нуля. Это предположение о симметрии $H(x)$ не обязательно, но позволяет удобно записать ответы в Теореме 2.2 далее.

Для полуинтервалов 
$$
B^+_j=(x_{j-1},x_j],\quad j=1,\ldots,m,\; m>1,\; 0=x_0<x_1<\ldots<x_m=\infty,
$$
пусть $
p^+_j:=G_0(x_{j})-G_0(x_{j-1})>0.\,\,$ Тогда при $H_{0}\,\Prob(\e_1 \in B^+_j)=p^+_j
$.\\
 Если ввести еще симметричные полуинтервалы 
 $ B^-_j=(-x_{j},-x_{j-1}]$, то\\
 $$
 \text{ при}\,\, H_{0}\,\,\Prob(\e_1 \in B^+_j \cup B^-_j)=2p^+_j:=p_j^0,\,p_j^0>0.\,\sum_{j=1}^m p_j^0=1.
 $$
Пусть $\hat \nu_j^+$ обозначает число остатков (они определены в (2.4)) среди $\{\hat \e_t,\,t=1,\ldots,n\}$, попавших в $B^+_j$, а 
$\hat \nu_j^-$ обозначает число остатков, попавших в $B^-_j$. Пусть
$$
\hat\nu_j=\hat \nu_j^+ +\hat \nu_j^-.
$$

Интересующая нас тестовая статистика типа хи-квадрат для $H_{0}$ имеет вид
$$
\hat{\chi}^2_{n} = \sum_{j=1}^m \frac{(\hat{\nu}_j - n p_j^0)^2}{n p_j^0}.
$$
Статистика $\hat{\chi}^2_{n} $ является простым функционалом от $\hat S_n(x)$. Действительно,

$$
n[\hat S_n(x_j)-\hat S_n(x_{j-1})]=\frac{1}{2}\{n[\hat G_n(x_j)-\hat G_n(x_{j-1})]+n[\hat G_n(-x_{j-1})-\hat G_n(-x_{j})]\}=
$$
$$
(\hat\nu^+_j +\hat\nu^-_j)/2)=\hat\nu_j/2.\eqno(2.10)
$$
Будем называть $\hat{\chi}^2_{n}$ симметризованной статистикой Пирсона.\\
Чтобы описать асимптотические свойства $\hat{\chi}^2_{n}$ нам понадобятся некоторые обозначения.
Введем диагональную матрицу
$$
\bm P_0=\mbox{diag}\{p_1^0,\ldots,p_m^0\},
$$
и вектора
$$
\bm p^0=(p_1^0,\ldots,p_m^0)^T;\quad
\bm p^H=(p_1^H,\ldots,p_m^H)^T\,\,\text{ c} \,\, p_j^H=2(H(x_j)-H(x_{j-1}));
$$
$$
\bm p^A=(p_1^A,\ldots,p_m^A)^T,\,\,p_j^A=A_n(x_j)-A_n(x_{j-1}).
$$
Нам понадобится еще вектор
$$
\bm \delta(\Pi)=(\delta_1(\Pi),\ldots,\delta_m(\Pi))^T\quad \text{с компонентами}\,\, \delta_j(\Pi)=2[\Delta_{S}(x_{j},\Pi)-\Delta_{S}(x_{j-1},\Pi)],
$$
сдвиг $\Delta_{S}(x,\Pi)$ определен в (2.8).

Пусть $ \nu_j^+$ обозначает число  $\{ \e_t,\,t=1,\ldots,n\}$, попавших в $B^+_j$, а 
$ \nu_j^-$  попавших в $B^-_j$. Пусть $\nu_j= \nu_j^+ + \nu_j^-$.

Из определения $S_n(x)$ и $\nu_j$ следует:
$$
n[ S_n(x_j)- S_n(x_{j-1})]=\nu_j/2.\eqno(2.11)
$$
В силу соотношений (2.10) --(2.11), Следствия 2.2 и определения вектора $\delta(\Pi)$ имеем при $H_{1n}(\rho)$:
$$
n^{1/2}(\frac{\hat\nu_j}{n}-
p_j^0)-n^{1/2}(\frac{\nu_j}{n}-p_j^0)=
n^{1/2}(\frac{\hat\nu_j}{n}-\frac{\nu_j}{n})=
$$
$$
2n^{1/2}\{[\hat S_n(x_j)-\hat S_n(x_{j-1}]-[S_n(x_j)-S_n(x_{j-1}]\}=\gamma \delta_j(\Pi)+o_P(1),\quad n\to \infty.\eqno(2.12)
$$
Здесь и далее $o_P(1)$ обозначает величину, сходящуюся к нулю по вероятности равномерно по $\gamma \leq \Gamma,\,\rho\leq R$.\\
Введем вектора 
$$
\bm\hat \nu=(\hat\nu_1,\ldots,\hat\nu_n)^T,\quad \bm \nu=(\nu_1,\ldots,\nu_n)^T.
$$
В силу (2.12) при $H_{1n}(\rho)$ имеем:
$$
n^{1/2}(\frac{\bm \hat \nu}{n}-\bm p^0)=
n^{1/2}(\frac{\bm \nu}{n}-\bm p^0)+\gamma \bm\delta(\Pi)+o_P(1)=
$$
$$
n^{1/2}(\frac{\bm \nu}{n}-\bm p^A)+n^{1/2}(\bm p^A-\bm p^0)+\gamma \bm\delta(\Pi)+o_P(1)=
$$
$$
n^{1/2}(\frac{\bm \nu}{n}-\bm p^A)+\rho (\bm p^H-\bm p^0)+\gamma \bm\delta(\Pi)+o_P(1),\quad n\to \infty.\eqno(2.13)
$$
Напомним еще известный факт относительно слабой сходимости вектора $\bm \nu$  (см. \cite{Bold.MMS19}, доказательство Теоремы 2.1):
$$
n^{1/2}(\frac{\bm \nu}{n}-\bm p^A) \to \bm N(\bm 0,\,\, \bm P_0-\bm p^0(\bm p^0)^T),\quad n\to \infty.\eqno(2.14)
$$
Слабая сходимость в (2.14) равномерна по $\rho\leq R$.

Соотношений (2.13) --(2.14) достаточно, чтобы следуя схеме доказательства Теоремы 2.3 в \cite{Bold.MMS19} доказать наше  основное утверждение -- Теорему 2.2. В ней $ F_{m-1}(x, \hat\lambda^2(\rho,\gamma, \Pi))$ означает ф.р. нецентрального хи-квадрата с $m-1$ степенями свободы и параметром нецентральности $\hat\lambda^2(\rho,\gamma, \Pi)$.
\begin{Th}
Пусть альтернатива  $H_{1n}(\rho)$ верна. Пусть выполнены Условия (i)--(ii) с симметричными $G_0$ и $H$. Тогда для любых конечных  неотрицательных $\Gamma,\,R$
$$
\sup_{x\in\mathbb{ R}^1, \gamma \le \Gamma,\rho\leq R} |\Prob(\hat{\chi}^2_{n} \le x) - F_{m-1}(x, \hat\lambda^2(\rho,\gamma, \Pi))| \to 0,\quad n \to \infty.
$$
Параметр нецентральности равен
$$
\hat\lambda^2(\rho,\gamma, \Pi) =
| \bm{ P}^{-1/2}_0[\rho(\bm p^H-\bm p^0)+\gamma \bm{\delta}(\Pi)]|^2.
$$
\end{Th}
В силу Теоремы 2.2 при $H_0$ и $\gamma=0$ предельное распределение нашей тестовой статистики $\hat{\chi}^2_{n}$ будет обычное (центральное) распределение хи-квадрат с $m-1$ степенями свободы.

Для заданного $0<\alpha<1$ ($\,\alpha$ --асимптотический уровень в схеме без засорений), мы будем отвергать $H_{0}$, когда
$$
\hat{\chi}^2_{n} > \chi_{m-1}(1-\alpha), \eqno(2.15)
$$
$\chi_{m-1}(1-\alpha)$ -- квантиль уровня $1-\alpha$ распределения хи-квадрат с $m-1$ степенью свободы. 
Мощность такого теста равна
$$
W_n(\rho,\gamma,\Pi)=P(\hat{\chi}^2_{n} > \chi_{m-1}(1-\alpha)).
$$
В силу Теоремы 2.2 эта мощность сходится при $n\to\infty$ равномерно по $\gamma\leq \Gamma,\,\rho\leq R$
к асимптотической мощности
$$
W(\rho,\gamma,\Pi)= 1 - F_{m-1}(\chi_{m-1}(1-\alpha), \hat\lambda^2(\rho,\gamma, \Pi)),\quad
 W(0,0,\Pi)= \alpha.
 $$
 Используя простое неравенство
 $$
 |F_k(x,\lambda_1^2)-F_k(x,\lambda_2^2)|\leq 2 \sup_{x \in \mathbb{ R}^1}\phi(x)|\lambda_1-\lambda_2|,
 $$
 $\phi(x)$ стандартная гауссовская плотность, и определение $\hat\lambda^2(\gamma, \Pi))$, получаем:
 $$
\sup_{\Pi,\rho\geq 0}|W(\rho,\gamma,\Pi)-W(\rho,0,\Pi)| \to 0,\quad \gamma \to 0.\eqno(2.16)
$$
Соотношение (2.16) означает асимптотическую качественную робастность теста (2.15). Такая робастность означает, что  $H_{0}$ для малых  $\gamma$ можно проверять примерно с асимптотическим уровнем $\alpha$, и тест будет иметь примерно такую же асимптотическую мощность при $H_{1n}(\rho)$ c $\rho>0$, как в схеме без засорений. И все это  независимо от распределения засорений $\Pi$.

\end{document}